# The False Discovery Rate for Statistical Pattern Recognition

### Clayton Scott[*] and Gowtham Bellala

*Electrical Engineering and Computer Science*
*University of Michigan, Ann Arbor*
*e-mail:* `{clayscot,gowtham}@umich.edu`

### Rebecca Willett[†]

*Electrical and Computer Engineering*
*Duke University*
*e-mail:* `willett@duke.edu`

**Abstract:** The false discovery rate (FDR) and false nondiscovery rate (FNDR) have received considerable attention in the literature on multiple testing. These performance measures are also appropriate for classification, and in this work we develop generalization error analyses for FDR and FNDR when learning a classifier from labeled training data. Unlike more conventional classification performance measures, the empirical FDR and FNDR are not binomial random variables but rather a ratio of binomials, which introduces challenges not addressed in conventional analyses. We develop distribution-free uniform deviation bounds and apply these to obtain finite sample bounds and strong universal consistency.

## Contents

[*]C. Scott was partially supported by NSF award no. 0830480.
[†]R. Willett was partially supported by DARPA grant no. HR0011-06-1-0026 and NSF CAREER award no. NSF-CCF-06-43947.

## 1. Introduction

When learning a classifier from labeled training data, minimizing the probability of misclassification is often unsatisfactory. In a variety of applications, such as screening medical images for cancerous lesions or detecting landmines, false positives and negatives have different impacts. False detections of targets are problematic because of the time, money, and other resources which are invariably wasted as a result. Missed detections, on the other hand, may result in loss of life or destruction. For this reason, a number of methods for cost-sensitive [3,12] and Neyman-Pearson [6,17,18] classification have been developed that allow the user to effect a tradeoff between false positive and negative rates.

The probability of error, false positive rate, and false negative rate are all performance measures that reflect the performance of a classifier on a *single* future test point. However, it is often the case that we desire to classify *multiple* future test points. In this situation, the false positive and negative rates may not be the most appropriate measures of performance. If a classifier has a false positive rate of say 5%, and 1000 negative test points (e.g., no target present) are observed, we expect 50 of them to be declared positive. This may be unacceptable, especially in situations where large costs are involved in investigating false positives.

This situation is similar to the multiple testing problem in hypothesis testing. Consequently, many of the ideas from multiple testing are applicable in the classification setting. The basic approach is to consider alternative measures of size and power that are better suited to simultaneous inference, and to design decision rules based on these new performance measures.

In this paper, we consider the false discovery rate (FDR) [20], which has emerged as the method of choice for quantifying error rates meaningfully in many multiple testing situations, with applications ranging from wavelet denoising [8] to neuroimaging [13] to the analysis of DNA microarrays [10]. Control of the FDR, i.e., the fraction of declared positives (discoveries) that are in fact negative, ensures that follow-up investigations into declared positives must return a certain yield of actual positives. Such control is vital in applications where follow-up studies are time or resource consuming.

Several researchers, spurred by the seminal work of [4], have studied FDR control in the context of multiple hypothesis testing by assuming known distributions of observed statistics under the null hypothesis. FDR control is then achieved, typically, by adjusting p-values through single step, step-up or step-down procedures. It is important to note that such procedures are not applicable in the statistical learning context because we do not assume knowledge of the null distribution and must instead rely upon training data.

We develop basic results on the analysis of generalization error in FDR controlled classification, including uniform deviation bounds, finite sample performance guarantees, and strong universal consistency. Unlike traditional perfor-mance probabilities, whose empirical versions are related to binomial random variables, empirical versions of FDR and FNDR are related to ratios of binomial variables. This necessitates the development of novel concentration inequalities and methods of analysis.

## 1.1. *Notation*

More formally, in this paper we consider the following scenario: Let $\mathcal{X}$ be a set and $Z = (X, Y)$ be a random variable taking values in $\mathcal{Z} = \mathcal{X} \times \{0, 1\}$. The variable $X$ corresponds to a pattern or feature vector and $Y$ to a class label associated with $X$; $Y = 0$ corresponds to the null hypothesis (e.g., that no target is present) and $Y = 1$ corresponds to the alternative hypothesis (e.g., that a target is present). The distribution on $Z$ is unknown and is denoted by $\mathbf{P}$. Assume we make $n$ independent and identically distributed training observations of $Z$, denoted $Z^n = (X_i, Y_i)_{i=1}^n$.

A classifier is a function $h : \mathcal{X} \longrightarrow \{0, 1\}$ mapping feature vectors to class labels. Let $\mathcal{H}$ denote a collection of different classifiers. A false discovery occurs when $h(X) = 1$ but the true label is $Y = 0$. Similarly, a false nondiscovery occurs when $h(X) = 0$ but $Y = 1$. We define the *false discovery rate* (FDR)

$$\mathcal{R}_D(h) := \begin{cases} \mathbf{P}(Y = 0 \,|\, h(X) = 1), & \text{if } \mathbf{P}(h(X) = 1) > 0, \\ \infty, & \text{else,} \end{cases}$$

and the *false nondiscovery rate* (FNDR)

$$\mathcal{R}_{ND}(h) := \begin{cases} \mathbf{P}(Y = 1 \,|\, h(X) = 0), & \text{if } \mathbf{P}(h(X) = 0) > 0, \\ \infty, & \text{else.} \end{cases}$$

## 1.2. *Related Concepts*

These definitions, which are natural in the classification setting, coincide with the so-called *positive* FDR/FNDR of Storey [21, 22], so named because it can be seen to equal the expected fraction of false discoveries/nondiscoveries, conditioned on a positive number of discoveries/nondiscoveries having been made. Storey makes some decision-theoretic connections to classification [22], but does not consider learning from data.

Storey's definition does not cover the case where the conditioning event has probability zero. We define FDR and FNDR in these cases to be infinity. Our convention has the effect of assigning high costs to classifiers that fail to make at least some discoveries (and nondiscoveries). This is consistent with the multiple testing perspective, where the goal is to generate interesting hypotheses for further examination. A classifier that makes no discoveries is of no use for such purposes. Further comments on the definition of FDR and FNDR are given after the proof of Theorem 2.

In certain communities, different terms embody the idea behind FDR. In the medical diagnostic testing literature, the *positive predictive value* (PPV) is defined as the "proportion of patients with positive test results who are correctly diagnosed" [1]. In database information retrieval problems, the *precision* is defined as the ratio of the number of relevant documents retrieved by a search to the total number of documents retrieved by a search [23]. Both PPV and precision are equal to 1 - FDR. Precision is discussed further is Section 5.2.

Finally, several researchers have recently investigated connections between multiple testing and statistical learning theory. McAllester's PAC-Bayesian learning theory may be viewed as an extension of multiple testing procedures to (possibly uncountably) infinite collections of hypotheses [16]. Blanchard and Fleuret present an extension of the Occam's razor principle for generalization error analysis in classification, and apply it to derive p-value adjustment procedures for controlling FDR [5]. Arlot et al. develop concentration inequalities that apply to multiple testing with correlated observations [2]. None of these works consider FDR/FNDR as performance criteria for classification.

### 1.3. Connections to Cost-Sensitive Learning

In Sections 3 and Section 4 we consider the performance measure $\mathcal{E}_\lambda(h) := \mathcal{R}_{ND}(h) + \lambda \mathcal{R}_D(h)$. It can be shown that the global minimizers of this criterion have the form

$$h(x) = \mathbf{1}_{\{\eta(x) \geq c\}} \tag{1}$$

for some $c$, where $\eta(x) := P(Y = 1 | X = x)$ and, if necessary, this family of classifiers is extended by a standard randomization argument if its receiver operating characteristic (ROC) is not concave. Storey [22] gives a proof for the case where the two class-conditional distributions are continuous. The classifiers in (1) are also the optimal classifiers for Bayes cost-sensitive learning. That is, they are also the minimizers of weighted Bayes costs of the form

$$\mathbf{P}(h(X) = 0, Y = 1) + \gamma \mathbf{P}(h(X) = 1, Y = 0),$$

$\gamma > 0$, where $c = 1/(1 + \gamma)$. Proof of this fact is a direct generalization of the case of the probability of error, when $\gamma = 1$ [7].

Unfortunately, existing analyses for cost-sensitive classification cannot be readily applied to our problem. Given $\lambda$, it is true that our criterion $\mathcal{E}_\lambda(h)$ can be minimized by performing cost-sensitive classification with a certain cost parameter $\gamma$. The *critical issue* is that $\gamma$ is an implicit function of $\lambda$, and cannot be determined a priori without knowledge of the underlying distribution. Therefore, when only data are given, applying existing cost-sensitive classification methods to our problem would require estimating $\gamma$. In practice, this would most likely entail learning cost-sensitive classifiers $\hat{h}_{\gamma_i}$ for some grid of values $\{\gamma_i\}$ that grows increasingly dense as $n \to \infty$. Then, the best of these candidates would be selected by minimizing an estimate of $\mathcal{E}_\lambda(h)$. Such a procedure would likely be expensive computationally. From an analytical standpoint, it seems plausible that generalization error analyses for cost-sensitive classification could be useful; however, the need to search for a $\gamma$ that approximately minimizes our criterion would significantly complicate the analysis. The objective of our work is to develop a much more direct approach, which does not require repeated cost-sensitive classification.

Therefore, the distinction between our problem and cost-sensitive classification is in some ways analogous to the difference between the Neyman-Pearson and Bayesian theories of hypothesis testing. Even though these two problems have a likelihood ratio as their optimal solution, the specific thresholds for the likelihood ratios are determined in very different ways depending on which criterion is employed. In our setting, the differences are further compounded by the fact that we are learning from data.

### *1.4. Overview*

In the next section we present and prove uniform deviation bounds for FDR and FNDR. In Section 3, we discuss performance measures based on FDR and FNDR, and in Section 4 we establish the strong universal consistency of a learning rule with respect to the measure $\mathcal{E}_\lambda$. Section 5 treats performance measures which constrain FDR, and the final section offers a concluding discussion. Several aspects of our analysis deviate from standard techniques, a consequence of certain unique features of FDR and FNDR, and we highlight these throughout the paper.

## 2. Uniform Deviation Bounds

Define empirical analogues to the FDR and FNDR according to

$$
\widehat{\mathcal{R}}_D(h) \quad := \quad \begin{cases} \frac{1}{n_D(h,Z^n)} \sum_{i=1}^n \mathbf{1}_{\{Y_i=0,h(X_i)=1\}}, & n_D(h,Z^n) > 0, \\ \infty, & n_D(h,Z^n) = 0, \end{cases}
$$

$$
\widehat{\mathcal{R}}_{ND}(h) \quad := \quad \begin{cases} \frac{1}{n_{ND}(h,Z^n)} \sum_{i=1}^n \mathbf{1}_{\{Y_i=1,h(X_i)=0\}}, & n_{ND}(h,Z^n) > 0, \\ \infty, & n_{ND}(h,Z^n) = 0, \end{cases}
$$

where $n_D(h,Z^n) = \sum_{i=1}^n \mathbf{1}_{\{h(X_i)=1\}}$ and $n_{ND}(h,Z^n) = \sum_{i=1}^n \mathbf{1}_{\{h(X_i)=0\}}$ are binomial random variables. This section describes a uniform bound on the amount by which the empirical estimate of FDR/FNDR can deviate from the true value. Note that unlike the usual empirical estimates for the probability of error/false positive rate/false negative rate, here both numerator and denominator are random, and both depend on $h$.

Assume $\mathcal{H}$ is countable, and let $[\![h]\!]$ be a real valued functional on $\mathcal{H}$ such that $\sum_{h \in \mathcal{H}} 2^{-[\![h]\!]} \leq 1$. Such a functional can be identified with a prefix code for $\mathcal{H}$, in which case $[\![h]\!]$ is the codelength associated to $h$. If $\sum_{h \in \mathcal{H}} 2^{-[\![h]\!]} = 1$, then $2^{-[\![h]\!]}$ may be viewed as a prior distribution on $\mathcal{H}$.

For $\delta > 0$, we introduce the *penalty* terms

$$
\phi_D(h, \delta) \quad := \quad
\begin{cases}
\sqrt{\dfrac{[\![h]\!]\log 2 + \log(2/\delta)}{2 n_D(h, Z^n)}}, & n_D(h, Z^n) > 0, \\
\infty, & n_D(h, Z^n) = 0,
\end{cases}
$$

$$
\phi_{ND}(h, \delta) \quad := \quad
\begin{cases}
\sqrt{\dfrac{[\![h]\!]\log 2 + \log(2/\delta)}{2 n_{ND}(h, Z^n)}}, & n_{ND}(h, Z^n) > 0, \\
\infty, & n_{ND}(h, Z^n) = 0.
\end{cases}
$$

The interpretation of these expressions as penalties comes from the learning algorithms studied below, where we minimize the empirical error plus a penalty to avoid overfitting. Note that the penalties are data dependent.

**Theorem 1.** *With probability at least $1 - \delta$ with respect to the draw of the training data,*

$$
|\mathcal{R}_D(h) - \widehat{\mathcal{R}}_D(h)| \leq \phi_D(h, \delta) \tag{2}
$$

*for all $h \in \mathcal{H}$. Similarly, with probability at least $1 - \delta$ with respect to the draw of the training data,*

$$
|\mathcal{R}_{ND}(h) - \widehat{\mathcal{R}}_{ND}(h)| \leq \phi_{ND}(h, \delta) \tag{3}
$$

*for all $h \in \mathcal{H}$. The results are independent of the underlying probability distribution.*

Because of the form of the penalty terms, the bound is larger for classifiers $h$ that are more complex, as represented through the codelength $[\![h]\!]$, and smaller when more discoveries/nondiscoveries are made. This result leads to finite sample bounds and strong universal consistency for certain learning rules based on minimization of the penalized empirical error, as developed in the sequel.

*Proof.* We prove the first statement, the second being similar. For added clarity, write the penalty as $\phi_D(h, \delta, n_D(h, Z^n))$, where

$$
\phi_D(h, \delta, k) :=
\begin{cases}
\sqrt{\dfrac{[\![h]\!]\log 2 + \log(2/\delta)}{2k}}, & k > 0, \\
\infty, & k = 0.
\end{cases}
$$

Consider a fixed $h \in \mathcal{H}$. The fundamental concentration inequality underlying our bounds is Hoeffding's [14], which, in one form, states that if $S_k$ is the sum of $k > 0$ independent random variables bounded between zero and one, and $\mu = \mathbf{E}[S_k]$, then

$$
\mathbf{P}(|\mu - S_k| > k\epsilon) \leq 2e^{-2k\epsilon^2}.
$$

To apply Hoeffding's inequality, we need the following conditioning argument. Let $V = (V_1, \ldots, V_n) \in \{0, 1\}^n$ be a binary indicator vector, with $V_i = \mathbf{1}_{\{h(X_i) = 1\}}$. Let $\mathcal{V}_k$ denote the set of all $v = (v_1, \ldots, v_n) \in \{0, 1\}^n$ such that $\sum_{i=1}^{n} v_i = k$. We may then write

$$\mathbf{P}(|\mathcal{R}_D(h) - \widehat{\mathcal{R}}_D(h)| > \phi_D(h, \delta, n_D(h, Z^n)))$$

$$= \sum_{k=0}^{n} \sum_{v \in \mathcal{V}_k} \mathbf{P}(|\mathcal{R}_D(h) - \widehat{\mathcal{R}}_D(h)| > \phi_D(h, \delta, k)|V = v)\mathbf{P}(V = v)$$

$$= \sum_{k=0}^{n} \sum_{v \in \mathcal{V}_k} \mathbf{P}(|k\mathcal{R}_D(h) - k\widehat{\mathcal{R}}_D(h)| > k\phi_D(h, \delta, k)|V = v)\mathbf{P}(V = v),$$

First note that $|\mathcal{R}_D(h) - \widehat{\mathcal{R}}_D(h)| \leq \phi_D(h, \delta)$ with probability one when $n_D(h, Z^n) = 0$. We now apply Hoeffding's inequality for each $k \geq 1$ and $v \in \mathcal{V}_k$, conditioning on $V = v$. Setting $S_k = k\widehat{\mathcal{R}}_D(h)$, we have

$$\begin{aligned}
\mu &= \mathbf{E}[S_k|V = v] \\
&= k\mathbf{E}[\widehat{\mathcal{R}}_D(h)|V = v] \\
&= \mathbf{E}[\sum_{i=1}^{n} \mathbf{1}_{\{Y_i=0, h(X_i)=1\}}|V = v] \\
&= \mathbf{E}[\sum_{i:v_i=1} \mathbf{1}_{\{Y_i=0\}}|V = v] \\
&= k\mathbf{P}(Y = 0|h(X) = 1) \\
&= k\mathcal{R}_D(h),
\end{aligned}$$

where in the next to last step we use independence of the realizations. Therefore, applying Hoeffding's inequality conditioned on $V = v \in \mathcal{V}_k$ yields

$$\mathbf{P}(|\mathcal{R}_D(h) - \widehat{\mathcal{R}}_D(h)| > \phi_D(h, \delta, n_D(h, Z^n)))$$

$$\leq \sum_{k=1}^{n} \sum_{v \in \mathcal{V}_k} 2e^{-2k\phi_D^2(h, \delta, k)}\mathbf{P}(V = v)$$

$$\leq \sum_{k=1}^{n} \sum_{v \in \mathcal{V}_k} \delta 2^{-\llbracket h \rrbracket}\mathbf{P}(V = v)$$

$$= \delta 2^{-\llbracket h \rrbracket}(1 - \mathbf{P}(\sum V_i = 0)) \leq \delta 2^{-\llbracket h \rrbracket}.$$

The result now follows by applying the union bound over all $h \in \mathcal{H}$. $\qquad\square$

The technique of conditioning on the random denominator of a ratio of binomials has also been applied in others settings [15, 18]. Unlike those works, however, here the binomial denominator depends on the classifier $h$. This presents difficulties for extending the above techniques to uncountable classes $\mathcal{H}$. See the final section for further discussion of this issue.

## 3. Measuring Performance

We would like to be able to make FDR/FNDR related guarantees about how a data-based classifier $\widehat{h}$ performs. For this, we need to specify a performance measure or optimality criterion that incorporates both FDR and FNDR quantities simultaneously. One possibility is to specify a number $0 < \alpha < 1$ and seek the classifier such that $\mathcal{R}_{ND}(h)$ is minimal while $\mathcal{R}_D(h) \leq \alpha$. We consider this setting in Section 5. Another is to specify a constant $\lambda > 0$ reflecting the relative cost of FDR to FNDR, and minimize

$$\mathcal{E}_\lambda(h) := \mathcal{R}_{ND}(h) + \lambda \mathcal{R}_D(h).$$

This measure was introduce by Storey [22], but was not studied in a learning context. The uniform deviation bounds of the previous section immediately imply the following computable bound on a classifier's performance with respect to this measure.

**Corollary 1.** *For any $\delta > 0$ and $n \geq 1$, with probability at least $1 - 2\delta$ with respect to the draw of the training data,*

$$\mathcal{E}_\lambda(h) \leq \widehat{\mathcal{R}}_{ND}(h) + \phi_{ND}(h,\delta) + \lambda[\widehat{\mathcal{R}}_D(h) + \phi_D(h,\delta)]$$

*for all $h \in \mathcal{H}$.*

In the next section, we analyze a learning rule based on minimizing the bound of Corollary 1, and establish its strong universal consistency.

## 4. Strong Universal Consistency

Denote the globally optimal value of the performance measure by

$$\mathcal{E}_\lambda^* := \inf_h \mathcal{E}_\lambda(h),$$

where the inf is over all measurable $h : \mathcal{X} \to \{0, 1\}$. We seek a learning rule $\widehat{h}_{\lambda,n}$ such that $\mathcal{E}_\lambda(\widehat{h}_{\lambda,n}) \to \mathcal{E}_\lambda^*$ almost surely, regardless of the underlying probability distribution. Thus let $\{\mathcal{H}_k\}_{k \geq 1}$ be a family of finite sets of classifiers with universal approximation capability. That is, assume that $\lim_{k \to \infty} \inf_{h \in \mathcal{H}_k} \mathcal{E}_\lambda(h) = \mathcal{E}_\lambda^*$ for all distributions on $(X, Y)$. Furthermore, assume this family to be *nested*, meaning $\mathcal{H}_1 \subseteq \mathcal{H}_2 \subseteq \mathcal{H}_3 \cdots$. For example, if $\mathcal{X} = [0, 1]^d$, we may take $\mathcal{H}_k$ to be the collection of histogram classifiers based on a binwidth of $2^{-k}$. Recall that we can set $[\![h]\!] = \log_2 |\mathcal{H}_k|$ for $h \in \mathcal{H}_k$, where $|\mathcal{H}_k|$ is the cardinality of $\mathcal{H}_k$. For histograms, we have $|\mathcal{H}_k| = 2^{2^{kd}}$ and hence $[\![h]\!] = 2^{kd} \log 2$.

The bound of Corollary 1 suggests bound minimization as a strategy for selecting a classifier empirically. However, rather than minimizing over all possible classifiers in some $\mathcal{H}_k$, we first discard those classifiers whose empirical numbers of discoveries or nondiscoveries are too small. In these cases, the penalties are possibly quite large, and we are unable to obtain tight concentrations of empirical FDR/FNDR measures around their true values. However, as $n$ increases, we are able to admit classifiers with increasingly small proportions of (non)discoveries, so that in the limit, we can still approximate arbitrary distributions. This aspect is another unique feature of FDR/FNDR compared to traditional performance measures.

Formally, set $\delta_n = 1/n^2$ and define

$$\widehat{\mathcal{H}}_n := \{h \in \mathcal{H}_{k_n} : \frac{n_{ND}(h, Z^n)}{n} \ge p_n, \frac{n_D(h, Z^n)}{n} \ge p_n\},$$

where $p_n := (\log n)^{-1}$. Here $k_n$ is such that $k_n \to \infty$ as $n \to \infty$ and $\log |\mathcal{H}_{k_n}| = o(n/\log n)$. For the histogram example, $\log |\mathcal{H}_{k_n}| = 2^{k_n d} \log 2$, and thus the assumed conditions on the growth of $k_n$ are essentially the same (up to a logarithmic factor) as for consistency of histograms in other problems. For example, in standard classification, $2^{k_n d} = o(n)$ is required [7].

Denote the bound of Corollary 1 by

$$\widehat{\mathcal{E}}_\lambda(h) := \widehat{\mathcal{R}}_{ND}(h) + \phi_{ND}(h, \delta_n) + \lambda[\widehat{\mathcal{R}}_D(h) + \phi_D(h, \delta_n)],$$

and define the classification rule

$$\widehat{h}_{\lambda,n} := \underset{h \in \widehat{\mathcal{H}}_n}{\arg\min} \ \widehat{\mathcal{E}}_\lambda(h).$$

If $\widehat{\mathcal{H}}_n = \emptyset$, then $\widehat{h}_{\lambda,n}$ may be defined arbitrarily.

**Theorem 2.** *For any distribution on $(X, Y)$, and any $\lambda > 0$,*

$$\mathcal{E}_\lambda(\widehat{h}_{\lambda,n}) \to \mathcal{E}_\lambda^*$$

*almost surely. That is, $\widehat{h}_{\lambda,n}$ is strongly universally consistent.*

*Proof.* First consider the case where there is no measurable $h : \mathcal{X} \to \{0, 1\}$ such that both $\mathbf{P}(h(X) = 0) > 0$ and $\mathbf{P}(h(X) = 1) > 0$. This occurs when $X$ is deterministic. Then $\mathcal{E}_\lambda^* = \infty$, and trivially $\widehat{h}_{\lambda,n}$ achieves optimal performance. So assume this is not the case.

By the Borel-Cantelli lemma [7, 9], it suffices to show that for each $\epsilon > 0$

$$\sum_{n=1}^\infty \mathbf{P}(\Omega^n) < \infty,$$

where

$$\Omega^n := \{Z^n : \mathcal{E}_\lambda(\widehat{h}_{\lambda,n}) - \mathcal{E}_\lambda^* \ge \epsilon\}.$$

Introduce the event

$$\Theta^n = \{Z^n : \widehat{\mathcal{H}}_n \ne \emptyset\}.$$

Then

$$\mathbf{P}(\Omega^n) = \mathbf{P}(\Omega^n | \Theta^n)\mathbf{P}(\Theta^n) + \mathbf{P}(\Omega^n | \overline{\Theta^n})\mathbf{P}(\overline{\Theta^n}) \le \mathbf{P}(\Omega^n | \Theta^n) + \mathbf{P}(\overline{\Theta^n}),$$

and therefore

$$\sum_{n=1}^\infty \mathbf{P}(\Omega^n) \le \sum_{n=1}^\infty \mathbf{P}(\Omega^n | \Theta^n) + \sum_{n=1}^\infty \mathbf{P}(\overline{\Theta^n}). \tag{4}$$

We will bound these two terms separately.

Consider the second term.

**Lemma 1.** *Let $\nu > 0$ and assume $\mathcal{E}_\lambda^* < \infty$. There exist $h'$ and $N_1$ such that $\mathcal{E}_\lambda(h') \leq \mathcal{E}_\lambda^* + \nu$ and, for all $n > N_1$, $\mathbf{P}(h' \in \widehat{\mathcal{H}}_n) \geq 1 - 1/n^2$.*

*Proof.* By the universal approximation assumption, there exists $m$ and $h' \in \mathcal{H}_{k_m}$ such that $\mathcal{E}_\lambda(h') \leq \mathcal{E}_\lambda^* + \nu$. Since $\mathcal{E}_\lambda^* < \infty$, this $h'$ necessarily has both $\mathbf{P}(h'(X) = 0) > 0$ and $\mathbf{P}(h'(X) = 1) > 0$. Denote $q := \min\{\mathbf{P}(h'(X) = 1), \mathbf{P}(h'(X) = 0)\} > 0$. Introduce

$$\tau_n := \sqrt{\frac{\log(2/\delta_n)}{2n}}.$$

By Hoeffding's inequality, with probability at least $1 - \delta_n$, $|\mathbf{P}(h'(X) = 1) - n_D(h', Z^n)/n| = |\mathbf{P}(h'(X) = 0) - n_{ND}(h', Z^n)/n| \leq \tau_n$. Since $\delta_n = 1/n^2$, we have that $\tau_n = o(p_n)$. Now choose $N_1$ such that $\tau_{N_1} \leq p_{N_1}$ and $2p_{N_1} \leq q$. Then, for a sample of size $n = N_1$, $\min\{n_D(h', Z^n)/N_1, n_{ND}(h', Z^n)/N_1\} \geq q - \epsilon_{N_1} \geq 2p_{N_1} - \tau_{N_1} \geq p_{N_1}$ with probability at least $1 - \delta_{N_1} = 1 - 1/N_1^2$. Since $p_n$ is decreasing and $\{\mathcal{H}_k\}$ is nested, the same is true for all $n > N_1$. $\square$

By this lemma we have $\mathbf{P}(\overline{\Theta^n}) \leq \delta_n = 1/n^2$ for all $n > N_1$ (Here we only the need the second part of the conclusion of the lemma; later we use the lemma in its full generality). Thus

$$\sum_{n=1}^\infty \mathbf{P}(\overline{\Theta^n}) \leq N_1 + \sum_{n > N_1} \frac{1}{n^2} < \infty.$$

Now consider the first term on the right-hand side of (4). Define the events

$$\Omega_1^n := \{Z^n : \mathcal{E}_\lambda(\widehat{h}_{\lambda,n}) - \inf_{h \in \widehat{\mathcal{H}}_n} \mathcal{E}_\lambda(h) \geq \frac{\epsilon}{2}\}$$

$$\Omega_2^n := \{Z^n : \inf_{h \in \widehat{\mathcal{H}}_n} \mathcal{E}_\lambda(h) - \mathcal{E}_\lambda^* \geq \frac{\epsilon}{2}\}$$

Since $\Omega^n \subset \Omega_1^n \bigcup \Omega_2^n$, we have

$$\sum_{n=1}^\infty \mathbf{P}(\Omega^n | \Theta^n) \leq \sum_{n=1}^\infty \mathbf{P}(\Omega_1^n | \Theta^n) + \sum_{n=1}^\infty \mathbf{P}(\Omega_2^n | \Theta^n). \tag{5}$$

We consider the two terms individually and show that each of them is finite.

To bound the first term on the right-hand side of (5) we use the following lemma.

**Lemma 2.** *If $\widehat{\mathcal{H}}_n \neq \emptyset$, then*

$$\mathcal{E}_\lambda(\widehat{h}_{\lambda,n}) - \inf_{h \in \widehat{\mathcal{H}}_n} \mathcal{E}_\lambda(h) \leq 2 \sup_{h \in \widehat{\mathcal{H}}_n} |\mathcal{E}_\lambda(h) - \widehat{\mathcal{E}}_\lambda(h)|.$$

*Proof.* Let $h' \in \widehat{\mathcal{H}}_n$ be arbitrary. By the definition of $\widehat{h}_{\lambda,n}$, $\widehat{\mathcal{E}}_\lambda(\widehat{h}_{\lambda,n}) \leq \widehat{\mathcal{E}}_\lambda(h')$. Hence

$$
\begin{aligned}
\mathcal{E}_\lambda(\widehat{h}_{\lambda,n}) &= \mathcal{E}_\lambda(\widehat{h}_{\lambda,n}) - \widehat{\mathcal{E}}_\lambda(\widehat{h}_{\lambda,n}) + \widehat{\mathcal{E}}_\lambda(\widehat{h}_{\lambda,n}) - \mathcal{E}_\lambda(h') + \mathcal{E}_\lambda(h') \\
&\leq \mathcal{E}_\lambda(\widehat{h}_{\lambda,n}) - \widehat{\mathcal{E}}_\lambda(\widehat{h}_{\lambda,n}) + \widehat{\mathcal{E}}_\lambda(h') - \mathcal{E}_\lambda(h') + \mathcal{E}_\lambda(h') \\
&\leq 2 \sup_{h \in \widehat{\mathcal{H}}_n} |\mathcal{E}_\lambda(h) - \widehat{\mathcal{E}}_\lambda(h)| + \mathcal{E}_\lambda(h').
\end{aligned}
$$

Since $h'$ was arbitrary, the result now follows. $\qquad\square$

Define the events

$$
\begin{aligned}
\Omega_{11}^n &:= \{Z^n : \sup_{h \in \widehat{\mathcal{H}}_n} |\mathcal{R}_{ND}(h) - \widehat{\mathcal{R}}_{ND}(h)| \geq \frac{\epsilon}{16}\} \\
\Omega_{12}^n &:= \{Z^n : \sup_{h \in \widehat{\mathcal{H}}_n} |\mathcal{R}_D(h) - \widehat{\mathcal{R}}_D(h)| \geq \frac{\epsilon}{16\lambda}\} \\
\Omega_{13}^n &:= \{Z^n : \sup_{h \in \widehat{\mathcal{H}}_n} |\phi_{ND}(h, \delta_n)| \geq \frac{\epsilon}{16}\} \\
\Omega_{14}^n &:= \{Z^n : \sup_{h \in \widehat{\mathcal{H}}_n} |\phi_D(h, \delta_n)| \geq \frac{\epsilon}{16\lambda}\}
\end{aligned}
$$

From Lemma 2 it follows that

$$
\Omega_1^n \subset \bigcup_{i=1}^4 \Omega_{1i}^n
$$

and hence it suffices to show

$$
\sum_{n=1}^\infty \mathbf{P}(\Omega_{1i}^n | \Theta^n)
$$

is finite for each $i = 1, 2, 3, 4$. We shall consider $\Omega_{11}$ and $\Omega_{13}$, the other two cases following similarly.

For $h \in \widehat{\mathcal{H}}_n$ we have $n_{ND}(h, Z^n)/n \geq p_n$ and therefore

$$
\begin{aligned}
\phi_{ND}(h, \delta_n) &= \sqrt{\frac{\log |\mathcal{H}_{k_n}| + \log(2n^2)}{2n_{ND}(h, Z^n)}} \\
&\leq \sqrt{(\log |\mathcal{H}_{k_n}| + \log(2n^2)) \frac{\log n}{2n}} < \frac{\epsilon}{16}
\end{aligned}
$$

for $n \geq N_2$, for some $N_2$ sufficiently large. Here we use $\delta_n = 1/n^2$ and $\log |\mathcal{H}_{k_n}| = o(n/\log n)$. Then

$$
\sum_{n=1}^\infty \mathbf{P}(\Omega_{13}^n | \Theta^n) \leq N_2.
$$

Furthermore, by the uniform deviation bound,

$$\sum_{n=1}^{\infty} \mathbf{P}(\Omega_{11}^n | \Theta^n) \leq N_2 + \sum_{n > N_2} \frac{1}{n^2} < \infty.$$

Now consider the event $\Omega_2^n$. Applying Lemma 1 with $\nu = \epsilon/2$, we have that

$$\sum_{n=1}^{\infty} \mathbf{P}(\Omega_2^n | \Theta_n) \leq N_1 + \sum_{n > N_1} \frac{1}{n^2} < \infty.$$

$\square$

In the definitions of $\mathcal{R}_D(h)$ and $\mathcal{R}_{ND}(h)$, we define these quantities to be infinity when the conditioning event has probability zero (see Introduction). This forces the globally optimal classifier to have both $\mathbf{P}(h(X) = 1) > 0$ and $\mathbf{P}(h(X) = 0) > 0$ whenever possible. The same property would hold provided $\mathcal{R}_D(h) > (1 + \lambda)/\lambda$ when $\mathbf{P}(h(X) = 1) = 0$ and $\mathcal{R}_{ND}(h) > (1 + \lambda)$ when $\mathbf{P}(h(X) = 0) = 0$. Were we to define FDR or FNDR to be smaller, our consistency argument would not apply universally. In particular, it might fail for distributions where the global minimizer of $\mathcal{E}_\lambda$ has either $\mathbf{P}(h(X) = 0) = 0$ or $\mathbf{P}(h(X) = 0) = 1$, such as when $X$ is deterministic. In a preliminary version of this work, we defined $\mathcal{R}_D(h)$ and $\mathcal{R}_{ND}(h)$ to be zero when the conditioning event is a null event, and were able to prove consistency under a very mild condition on the underlying distribution [19].

## 5. Constraining FDR

In this section we apply Theorem 1 to analyze a rule that seeks to minimize the FNDR subject to the constraint that FDR $\leq \alpha$, where $\alpha$ is a user-defined significance level. In fact, we first present a more general result, and then deduce results for this and other constrained learning problems as corollaries.

Thus, let $\mathcal{H}$ be a collection of classifiers as before, but not necessarily finite. Let $\mathcal{R}_0$ and $\mathcal{R}_1$ be measures of Type I and Type II error. For example, these may be FDR and FNDR, false positive rate and false negative rate, or some combination thereof. Assume that for $i = 0, 1$, there exists a data-based estimate $\widehat{\mathcal{R}}_i$ of $\mathcal{R}_i$, and a penalty $\phi_i(h, \delta)$, which define a symmetric confidence interval for $\mathcal{R}_i$. That is, suppose that for any $0 < \delta < 1$,

$$\mathbf{P}_{Z^n}(\sup_{h \in \mathcal{H}} [|\mathcal{R}_i(h) - \widehat{\mathcal{R}}_i(h)| - \phi_i(h, \delta)] > 0) \leq \delta.$$

For $0 < \alpha < 1$ define

$$\begin{aligned} h_{\mathcal{H}, \alpha}^* &= \underset{h \in \mathcal{H}}{\arg\min}\ \mathcal{R}_1(h) \\ &\quad \text{s. t. } \mathcal{R}_0(h) \leq \alpha. \end{aligned}$$

Consider the learning rule

$$\widehat{h}_{\mathcal{H},\alpha} \;\;=\;\; \underset{h \in \mathcal{H}}{\arg\min}\;\; \widehat{\mathcal{R}}_1(h) + \phi_1(h,\delta) \qquad (6)$$

$$\text{s. t. } \widehat{\mathcal{R}}_0(h) \leq \alpha + \phi_0(h,\delta).$$

**Theorem 3.** *The learning rule defined in Eqn. (6) is such that, for any $\delta > 0$ and any $n \geq 1$, with probability at least $1 - 2\delta$ with respect to the draw of the training data,*

$$\mathcal{R}_1(\widehat{h}_{\mathcal{H},\alpha}) \leq \mathcal{R}_1(h^*_{\mathcal{H},\alpha}) + 2\phi_1(h^*_{\mathcal{H},\alpha},\delta)$$

*and*

$$\mathcal{R}_0(\widehat{h}_{\mathcal{H},\alpha}) \leq \alpha + 2\phi_0(\widehat{h}_{\mathcal{H},\alpha},\delta).$$

The result holds regardless of the data-generating distribution.

*Proof.* Assume that both

$$|\mathcal{R}_0(h) - \widehat{\mathcal{R}}_0(h)| \leq \phi_0(h,\delta) \quad \text{for all } h \in \mathcal{H} \qquad (7)$$

and

$$|\mathcal{R}_1(h) - \widehat{\mathcal{R}}_1(h)| \leq \phi_1(h,\delta) \quad \text{for all } h \in \mathcal{H}, \qquad (8)$$

which, by assumption, occurs with probability at least $1 - 2\delta$. By (7), we deduce the second half of the theorem from

$$\mathcal{R}_0(\widehat{h}_{\mathcal{H},\alpha}) \leq \widehat{\mathcal{R}}_0(\widehat{h}_{\mathcal{H},\alpha}) + \phi_0(\widehat{h}_{\mathcal{H},\alpha},\delta) \leq \alpha + 2\phi_0(\widehat{h}_{\mathcal{H},\alpha},\delta),$$

where the second inequality follows from $\widehat{\mathcal{R}}_0(\widehat{h}_{\mathcal{H},\alpha}) \leq \alpha + \phi_0(\widehat{h}_{\mathcal{H},\alpha},\delta)$, which follows from the definition of $\widehat{h}_{\mathcal{H},\alpha}$. To get the first half of the theorem, observe that $\widehat{\mathcal{R}}_0(h^*_{\mathcal{H},\alpha}) \leq \mathcal{R}_0(h^*_{\mathcal{H},\alpha}) + \phi_0(h^*_{\mathcal{H},\alpha},\delta) \leq \alpha + \phi_0(h^*_{\mathcal{H},\alpha},\delta)$. Therefore, $h^*_{\mathcal{H},\alpha}$ is among the candidates in the minimization defining $\widehat{h}_{\mathcal{H},\alpha}$. Then

$$
\begin{aligned}
\mathcal{R}_1(\widehat{h}_{\mathcal{H},\alpha}) &\leq& \widehat{\mathcal{R}}_1(\widehat{h}_{\mathcal{H},\alpha}) + \phi_1(\widehat{h}_{\mathcal{H},\alpha},\delta) \\
&\leq& \widehat{\mathcal{R}}_1(h^*_{\mathcal{H},\alpha}) + \phi_1(h^*_{\mathcal{H},\alpha},\delta) \\
&\leq& \mathcal{R}_1(h^*_{\mathcal{H},\alpha}) + 2\phi_1(h^*_{\mathcal{H},\alpha},\delta).
\end{aligned}
$$

$\square$

This theorem can immediately be combined with Theorem 1 to give performance guarantees for the case $\mathcal{R}_0(h) = \mathcal{R}_D(h)$ and $\mathcal{R}_1(h) = \mathcal{R}_{ND}(h)$, for a countable class $\mathcal{H}$. In particular, define the rule

$$\widehat{h}_{\mathcal{H},\alpha} \;\;=\;\; \underset{h \in \mathcal{H}}{\arg\min}\;\; \widehat{\mathcal{R}}_{ND}(h) + \phi_{ND}(h,\delta) \qquad (9)$$

$$\text{s. t. } \widehat{\mathcal{R}}_D(h) \leq \alpha + \phi_D(h,\delta).$$

We have the following.

**Corollary 2.** *Assume $\mathcal{H}$ is countable. For any $\delta > 0$ and any $n \geq 1$, with probability at least $1 - 2\delta$ with respect to the draw of the training data, the learning rule in (9) satisfies*

$$\mathcal{R}_{ND}(\widehat{h}_{\mathcal{H},\alpha}) \leq \mathcal{R}_{ND}(h^*_{\mathcal{H},\alpha}) + 2\phi_{ND}(h^*_{\mathcal{H},\alpha}, \delta)$$

*and*

$$\mathcal{R}_D(\widehat{h}_{\mathcal{H},\alpha}) \leq \alpha + 2\phi_D(\widehat{h}_{\mathcal{H},\alpha}, \delta).$$

To extend such a result to a universally consistent estimator, based on the discussion of Theorem 2, it would be necessary to take $\mathcal{H}$ growing with the sample size $n$, and to exclude classifiers making too few discoveries or nondiscoveries. The details are similar to those of Section 4, and a formal development is omitted.

### 5.1. Neyman-Pearson Classification

If we take $\mathcal{R}_0$ and $\mathcal{R}_1$ to be the false positive rate and false negative rate, respectively, we may apply Theorem 3 to recover and generalize known results for Neyman-Pearson classification [6, 18]. Specifically, set

$$\begin{aligned}
\mathcal{R}_{FP}(h) &:= \mathbf{P}(h(X) = 1 \,|\, Y = 0) \\
\mathcal{R}_{FN}(h) &:= \mathbf{P}(h(X) = 0 \,|\, Y = 1).
\end{aligned}$$

There are several possible penalties that provide uniform bounds on the deviation between these quantities and their natural empirical estimates,

$$\widehat{\mathcal{R}}_{FP}(h) := \begin{cases} \frac{1}{n_0} \sum_{i=1}^n \mathbf{1}_{\{Y_i = 0, h(X_i) = 1\}}, & n_0 > 0, \\ 0, & n_0 = 0, \end{cases}$$

$$\widehat{\mathcal{R}}_{FN}(h) := \begin{cases} \frac{1}{n_1} \sum_{i=1}^n \mathbf{1}_{\{Y_i = 1, h(X_i) = 0\}}, & n_1 > 0, \\ 0, & n_1 = 0, \end{cases}$$

where $n_j := \sum_{i=1}^n \mathbf{1}_{\{Y_i = j\}}$. Examples of such penalties (e.g., VC and Rademacher penalties) are given in [17]. As a concrete example, we state a result here for the case of countable $\mathcal{H}$. Thus define the penalties

$$\phi_{FP}(h, \delta) = \begin{cases} \sqrt{\frac{[\![h]\!] \log 2 + \log(2/\delta)}{2n_0}}, & n_0 > 0, \\ 1, & n_0 = 0, \end{cases}$$

$$\phi_{FN}(h, \delta) = \begin{cases} \sqrt{\frac{[\![h]\!] \log 2 + \log(2/\delta)}{2n_1}}, & n_1 > 0, \\ 1, & n_1 = 0. \end{cases}$$

Define the rule

$$\begin{aligned}
\widehat{h}_{\mathcal{H},\alpha} &= \operatorname*{arg\,min}_{h \in \mathcal{H}} \; \widehat{\mathcal{R}}_{FN}(h) + \phi_{FN}(h, \delta) \qquad\qquad (10) \\
&\quad \text{s. t. } \widehat{\mathcal{R}}_{FP}(h) \leq \alpha + \phi_{FP}(h, \delta).
\end{aligned}$$

We have the following.

**Corollary 3.** *Assume $\mathcal{H}$ is countable. For any $\delta > 0$ and any $n \geq 1$, with probability at least $1 - 2\delta$ with respect to the draw of the training data, the learning rule in (10) satisfies*

$$\mathcal{R}_{FN}(\widehat{h}_{\mathcal{H},\alpha}) \leq \mathcal{R}_{FN}(h^*_{\mathcal{H},\alpha}) + 2\phi_{FN}(h^*_{\mathcal{H},\alpha}, \delta)$$

*and*

$$\mathcal{R}_{FP}(\widehat{h}_{\mathcal{H},\alpha}) \leq \alpha + 2\phi_{FP}(\widehat{h}_{\mathcal{H},\alpha}, \delta).$$

We note that Theorem 3, applied in the context of Neyman-Pearson classification, is a stronger result than those in [6, 18], which do not explicitly allow penalties that depend on the classifier $h$.

### 5.2. Precision and Recall

As a final application of Theorem 3, we analyze the precision and recall performance measures, common in database information retrieval problems (see Introduction). Precision and recall can both be defined in terms of quantities already discussed. Denote the precision

$$\mathcal{Q}_{PR}(h) := \mathbf{P}(Y = 1 \,|\, h(X) = 1) := 1 - \mathcal{R}_D(h)$$

and the recall

$$\mathcal{Q}_{RE}(h) := \mathbf{P}(h(X) = 1 \,|\, Y = 1) = 1 - \mathcal{R}_{FN}(h),$$

and let $\widehat{\mathcal{Q}}_{PR}(h) := 1 - \widehat{\mathcal{R}}_D(h)$ and $\widehat{\mathcal{Q}}_{RE}(h) := 1 - \widehat{\mathcal{R}}_{FN}(h)$ be the empirical estimates. In this setting the goal is to find the classifier with the largest precision, while maintaining a recall of at least $\beta$, where $\beta$ is a user-specified level. Thus the optimal classifier in a given class $\mathcal{H}$ is

$$
\begin{aligned}
h^*_{\mathcal{H},\beta} &= \underset{h \in \mathcal{H}}{\arg\max} \ \mathcal{Q}_{PR}(h) \\
&\text{s. t. } \mathcal{Q}_{RE}(h) \geq \beta.
\end{aligned}
$$

Define the rule

$$
\begin{aligned}
\widehat{h}_{\mathcal{H},\beta} &= \underset{h \in \mathcal{H}}{\arg\max} \ \widehat{\mathcal{Q}}_{PR}(h) - \phi_D(h, \delta) \\
&\text{s. t. } \widehat{\mathcal{Q}}_{RE}(h) \geq \beta - \phi_{FN}(h, \delta).
\end{aligned}
\tag{11}
$$

We have the following.

**Corollary 4.** *Assume $\mathcal{H}$ is countable. For any $\delta > 0$ and any $n \geq 1$, with probability at least $1 - 2\delta$ with respect to the draw of the training data, the learning rule in (11) satisfies*

$$\mathcal{Q}_{PR}(\widehat{h}_{\mathcal{H},\beta}) \geq \mathcal{Q}_{PR}(h^*_{\mathcal{H},\beta}) - 2\phi_D(h^*_{\mathcal{H},\beta}, \delta)$$

*and*

$$\mathcal{Q}_{RE}(\widehat{h}_{\mathcal{H},\beta}) \geq \beta - 2\phi_{FN}(\widehat{h}_{\mathcal{H},\beta}, \delta).$$

*Proof.* To apply Theorem 3, note that maximizing $\mathcal{Q}_{PR}(h)$ is equivalent to minimizing $\mathcal{R}_D(h)$, that the constraint $\mathcal{Q}_{RE}(h) \geq \beta$ is equivalent to $\mathcal{R}_{FN}(h) \leq \alpha := 1 - \beta$, and similarly for the empirical objective and constraint. Furthermore, since $|\mathcal{Q}_{PR}(h) - \widehat{\mathcal{Q}}_{PR}(h)| = |\mathcal{R}_D(h) - \widehat{\mathcal{R}}_D(h)|$, and $|\mathcal{Q}_{RE}(h) - \widehat{\mathcal{Q}}_{RE}(h)| = |\mathcal{R}_{FN}(h) - \widehat{\mathcal{R}}_{FN}(h)|$, we have that the assumptions of Theorem 3 are satisfied with the stated penalties. □

## 6. Conclusion

This paper demonstrates that FDR and FNDR control is possible in the context of statistical learning theory, where the distribution of $(X, Y)$ is unknown except through training data. We develop empirical estimates of these quantities and derive uniform deviation bounds which assess the closeness of these empirical estimates to the true FDR and FNDR. Unlike most other performance measures in statistical learning theory, which are related to binomial random variables, the FDR and FNDR measures are related to ratios of binomial random variables, which requires the development of novel bounding techniques. These bounds are then used to analyze learning rules that minimize a weighted combination of FDR and FNDR, or that minimize FNDR subject to a constraint on FDR. Our strong universal consistency result indicates that it is necessary to prevent the learning algorithm from selecting classifiers making too few discoveries or nondiscoveries, as error estimates for such classifiers may be poor.

Extending our results to uncountable classes $\mathcal{H}$ is an interesting open question, and may require the development of new techniques. The standard proofs of common generalization error bounds for uncountable classes, such as Rademacher and VC penalties, rely on the introduction of an artificial "ghost" sample [7]. That technique would require every $h \in \mathcal{H}$ to have the same empirical number of discoveries (or nondiscoveries) on both the original and ghost samples, which is generally not the case. Recently El-Yaniv and Pechyony [11] have extended the ghost sample technique to cases where the training and ghost samples have different sizes (their results are stated in the context of transductive learning), and some of their arguments may be useful in this regard.

## References

[1] ALTMAN, D. G. AND BLAND, J. M. (1994). Diagnostic tests 2: predictive values. *Brit. Med. J. 309*, 102.

[2] ARLOT, S., BLANCHARD, G., AND ROQUAIN, E. (2007). Resampling-based confidence regions and multiple tests for a correlated random vector. In *Learning Theory: 20th Annual Conference on Learning Theory, COLT 2007*, N. H. B. C. Gentile, Ed. Springer-Verlag, Heidelberg, 127–141.

[3] BACH, F. R., HECKERMAN, D., AND HORVITZ, E. (2006). Considering cost asymmetry in learning classifiers. *J. Machine Learning Research*, 1713–1741.

[4] BENJAMINI, Y. AND HOCHBERG, Y. (1995). Controlling the false discovery rate: a practical and powerful approach to multiple testing. *J. R. Statist. Soc B* **57**, 1, 289–300.

[5] BLANCHARD, G. AND FLEURET, F. (2007). Occam's hammer. In *Learning Theory: 20th Annual Conference on Learning Theory, COLT 2007*, N. H. B. C. Gentile, Ed. Springer-Verlag, Heidelberg, 112–126.

[6] CANNON, A., HOWSE, J., HUSH, D., AND SCOVEL, C. (2002). Learning with the Neyman-Pearson and min-max criteria. Tech. Rep. LA-UR 02-2951, Los Alamos National Laboratory. http://www.c3.lanl.gov/ml/pubs/2002_minmax/paper.pdf.

[7] DEVROYE, L., GYÖRFI, L., AND LUGOSI, G. (1996). *A Probabilistic Theory of Pattern Recognition*. Springer, New York.

[8] DONOHO, D. AND JIN, J. (2004). Higher criticism for detecting sparse heterogeneous mixtures. *Ann. Stat.* **32**, 3, 962–994.

[9] DURRETT, R. (1991). *Probability: Theory and Examples*. Wadsworth & Brooks/Cole, Pacific Grove, CA.

[10] EFRON, B., TIBSHIRANI, R., STOREY, J., AND TUSHER, V. (2001). Empirical Bayes analysis of a microarray experiment. *Journal of the American Statistical Association 96*, 1151–1160.

[11] EL-YANIV, R. AND PECHYONY, D. (2007). Transductive Rademacher complexity and its applications. In *Proc. 20th Annual Conference on Learning Theory, COLT 2007*, N. Bshouty and C. Gentile, Eds. Springer-Verlag, Heidelberg, 157–171.

[12] ELKAN, C. (2001). The foundations of cost-sensitive learning. In *Proceedings of the 17th International Joint Conference on Artificial Intelligence*. Seattle, Washington, USA, 973–978.

[13] GENOVESE, C. R., LAZAR, N. A., AND NICHOLS, T. E. (2002). Thresholding of statistical maps in functional neuroimaging using the false discovery rate. *NeuroImage 15*, 870–878.

[14] HOEFFDING, W. (1963). Probability inequalities for sums of bounded random variables. *J. Amer. Statist. Assoc. 58*, 13–30.

[15] MANSOUR, Y. AND MCALLESTER, D. (2000). Generalization bounds for decision trees. In *Proceedings of the Thirteenth Annual Conference on Computational Learning Theory*, N. Cesa-Bianchi and S. Goldman, Eds. Palo Alto, CA, 69–74.

[16] MCALLESTER, D. (1999). Some PAC-Bayesian theorems. *Machine Learning 37*, 3, 355–363.

[17] SCOTT, C. (2007). Performance measures for Neyman-Pearson classification. *IEEE Trans. Inform. Theory 53*, 8, 2852–2863.

[18] SCOTT, C. AND NOWAK, R. (2005). A Neyman-Pearson approach to statistical learning. *IEEE Trans. Inform. Theory 51*, 8, 3806–3819.

[19] SCOTT, C. D., BELLALA, G., AND WILLETT, R. (2007). Generalization error analysis for FDR controlled classification. In *Proc. IEEE Workshop on Statistical Signal Processing*. Madison, WI.

[20] SORIC, B. (1989). Statistical discoveries and effect-size estimation. *J. Amer. Statist. Assoc. 84*, 608–610.

[21] STOREY, J. (2002). A direct approach to false discovery rates. *Journal of the Royal Statistical Society, Series B 64*, 479–498.

[22] STOREY, J. (2003). The positive false discovery rate: A Bayesian interpretation of the $q$-value. *Annals of Statistics 31:6*, 2013–2035.

[23] VAN RIJSBERGEN, C. J. (1979). *Information Retrieval*, 2nd ed. Butterworths, London.